# *Symbolic Languages and Ars Combinatoria*

*Godofredo Iommi Amunátegui*
Instituto de Física
Pontificia Universidad Católica de Valparaíso
Casilla 4059, Valparaíso-Chile
giommi@ucv.cl

Abstract

This article analyses some paragraphs of the *Dissertatio de Arte Combinatoria* (1666) where G.W. Leibniz considers the syntax of a language with a given number of primitive terms. We propose a new formulation which generalizes the philosopher's conception of such a formal system.

Kew-words: Languaje, symbol, combinatorics, syntax, semantics.

# Introduction

In a passage of the many-sided masterpiece composed in his youth, the *Dissertatio de Arte Combinatoria* (1666)[1], Leibniz considers a symbolic system and formulates, with clairvoyance, albeit indirectly, a characteristic of formal languages, i.e., the bi-univocal relation which may be established between a definite signification and a multiplicity of modes of expression. This fact appears as a remarkable achievement of the *Ars Combinatoria* viewed as "a general science that teaches a merely syntactical manipulation of signs"[2].

In the first part of this note, we present the text and we add a commentary where Leibniz's exposition is subjected to scrutiny. Some inconsistencies blur, so to say, Leibniz's insight on the structure of his own language. In the second section, we propose a new formulation which, in our opinion, clarifies these aspects of the main idea. We conclude the article with a brief consideration.

At this place a concise glossary is useful. Leibniz uses the term *complexiones* (complexions) to designate the combinations. When all possible combinations, without repetitions, are meant he uses the expression *complexiones simpliciter*. When referring to a combination of a special class, he writes *com2natio* (combination), *com3natio, con4natio*, etc. The one element complexions are called *1-nions*. The number of *complexiones simpliciter* of $k$ elements is $2^k-1$. The primitive terms are the elements of a class and may be designated by numbers. So (1, 2, 3) represents a class of three primitive terms. The derived terms are made up of primitive terms: (1,3) is a derived term. The number of primitive terms of a derived term corresponds to the exponent of the complexion and determines the distance of the derived term. The exponent and the distance, as well, of (1,3) is two.

---

[1] G.W. Leibniz, Dissertatio de Arte Combinatoria, Die Philosophischen Schriften, Ed. Gerhardt, IV, Georg Olms, Hildesheim, 1978; pp. 27-104.
[2] Eberhard Knobloch, Leibniz between *ars characteristica* and *ars inveniendi*: Unknown news about Cajori's "master-builder of mathematical notations", in Philosophical Aspects of Symbolic Reasoning in Early Modern Mathematics, Ed. A. Heeffer and M. Van Dyck, Ed. Dov Gabbay, Vol. 26, Ghent, Belgium, 2010; pp. 289-302.

# I
## A.- The text

Hereafter, we consider the pertinent paragraphs of the *Dissertatio*[3]. Our translation is by no means literal and is intended as a rough guide to the original:

65.2   All the primitive terms considered are in a class and they are designated by some notes; the most convenient manner will be to number them.

66.3   The primitive terms comprise things and, as well, modes or relations[4].

67.4   All the derived terms are at a different distance from the primitive ones, depending on how many terms compose them, that is depending on the exponent of the complexion, hence they must form the same number of classes as the exponents, and in the same class must be assembled the terms composed of the same number of primitive terms.

68.5   The terms derived by com2nation can only be written by writing the primitive terms which compose them; and, as the primate terms have been designated by numbers, the two numbers which indicate the two terms are written.

69.6   The terms derived by con3nation or even by other complexions having a greater exponent, i.e., terms which are in the third class and the following ones, may be written, each one, in as many different ways as the number of *complexiones simpliciter* corresponding to the exponent of those very same terms, the exponent-this time-is considered as the number of things. All this is based on the Application IX; as an example: let these primitive terms be noted by the numbers 3.6.7.9; and let the derived term of the third class be composed by con3nation, that is by three simple terms: 3.6.9, and let the second class be these com2nations: (1)3.6 (2)3.7 (3)3.9 (4)6.7 (5)6.9 (6)7.9; I say that such a term of the third class may be written as: 3.6.9, by means of all the simple terms; or by one simple term, and instead of the other two simple ones, writing the com2nation in this way: $\frac{1}{2}^9$, or as $\frac{3}{2}^6$, or as $\frac{5}{2}^3$. Afterwards it will be said what do these kind of fractions mean. The greater the distance of a class from the first one, greater the variation will be. For the terms of a previous class behave like subaltern *genera* in respect to some terms of the next variation.

---

[3] G.W. Leibniz, Op. Cit., p. 65; Appendix I.
[4] Leibniz indicates the distinction between concrete terms ("flower", etc.) and abstract terms ("quantity", etc.). See Massimo Mugnai, Leibniz e la logica simbolica, Sansoni, Firenze, 1973; p. 29.

70.7 The number of times that a derived term is cited out of its class is written as a fraction, where the upper number, that is, the numerator is the place number in the class and the lower term or denominator is the class number 8. Concerning the derived terms it is more convenient to write the intermediate terms instead of all the primitive ones, and through a careful consideration, from these the other ones will appear. However, the most fundamental step is to write all the primitive ones.

## B.- Commentary

Leibniz states that[5]:

(i) The terms derived by com2nation are indicated by two primitive terms, for example, $a_1$ and $a_2$.

(ii) The terms derived by con3nation are represented according to an instance[6] of the Application IX which deals with genus and species: The case of the genus ($a_1$, $a_2$, $a_3$) where $a_1$, $a_2$ and $a_3$ are the species. The genus is a con3nation, the species are 1-nions and the com2nations correspond to "intermediate genera". There are seven terms (1 con3nation + 3com2nations +3unions). Hence the total number of terms is $2^3-1$ (*complexiones simpliciter*).

Moreover,

(iii) To denote the derived terms of a con4nation, Leibniz introduces a system of symbols that juxtapose integers and fractions[7]. He considers four primitive terms ($a_1$, $a_2$, $a_3$, $a_4$). The derived terms of the third class (con3nations) are ($a_1$, a₂, a₃), (a₁, a₂, a₄), (a₁, a₃, a₄), (a₂, a₃, a₄), and the derived terms of the second class (com2nations) are: (1)a₁a₂, (2)a₁a₃, (3)a₁a₄, (4)a₂a₃, (5)a₂a₄, (6)a₃a₄. So a₁a₂ is the first term of the second class, a₁a₃ is the second term of the second class and so on. Any com2nation can be written as a fraction whose denominator is the class number and whose numerator is the number assigned to the com2nation in the class, Hence a₁a₄ corresponds to the fraction $\frac{3}{2}$. A con3nation is denoted by means of the juxtaposition of a fraction and an integer: $\frac{1}{2}a_3$ represents a₁a₂a₃. As a result (a₁, a₂, a₃); $\frac{1}{2}a_3$, $\frac{2}{2}a_2$ and $\frac{4}{2}a_1$ are equivalent symbolic representations.

---
[5] We slightly modify the original notation.
[6] G.W. Leibniz, Op. Cit., p. 60.
[7] We are grateful to Eberhard Knobloch for many illuminating conversations concerning this case.

Following this procedure 16 expressions can be derived from the four primitive terms (a₁, a₂, a₃, a₄).

(iv) From (i), (ii) and (iii), we may deduce that (a) for a com2nation the number of derived terms is 2. (b) for a con3nation the number of derived terms is 7. (c) for a con4nation the number of derived terms is 16.

Leibniz affirms that for every con*k*nation, the number of derived terms is $2^k-1$, i.e., the number of *complexiones simpliciter*. But such statement is valid only for *k*=3. A consistent theoretical ground seams to the lacking to carry out the language's construction.

## II
## Leibniz's idea revisited: a new formulation

Leibniz envisages a language with *k* primitive elements and $2^k-1$ derived terms. Let us examine anew a section of a paragraph[8] which gives a clue to solve the puzzle. Therein Leibniz exposes this procedure: A con4nation is decomposed into com2nations. The con4nation is written in terms of con3nations which are, in turn, represented as quasi-fractions. As a result, the 4 primitive terms are the building blocks of a language having 16 terms.

---
[8] Paragraph 69.6 (see part I(A)).

With the aim of achieving our purpose, we shall generalize Leibniz's method. For $k$ primitive terms we have (1) a con$k$nation (2) $C_{k-2}^{k}$ con($k$-2)nations (3) $C_{k-1}^{k}$ con($k$-1)nations. These con($k$-1)nations are expressed in semi-fractional notation by means of the con($k$-2)nations. There are $C_{k-1}^{k} \cdot C_{k-2}^{k-1}$ semi-fractional terms.

Proposition: The $k$ primitive terms are the basic elements of a language of $k^2$ symbols.

The proof of this statement is direct. From (3) we may deduce that the total number of symbols of the language is:

$$C_{k-1}^{k} + C_{k-1}^{k} C_{k-2}^{k-1} = C_{k-1}^{k}(1 + C_{k-2}^{k-1})$$

$$C_{k-1}^{k}(1 + C_{k-2}^{k-1}) = kC_{k-1}^{k} = k^2$$

Let us consider the case $k = 2$, i.e., $(a_1, a_2)$. We number the 1-nions: (1)$a_1$; (2)$a_2$. Hence the resulting language has 4 signs: $C_1^2 + C_1^2 \cdot C_0^1$

| $a_1$ | $\frac{1}{1}a_2$ |
|---|---|
| $a_2$ | $\frac{2}{1}a_1$ |

We have chosen this example to point out the presence of the combinatiorial number $C_0^k$ that may appear surprising in such a mathematical context. Nonetheless, Leibniz mentions it in the *Dissertatio*[9]. Knobloch and Serres already noted its presence[10]. We develop the case $k$=6 to exemplify a more complex language. We have:
(i) 15 con4nations: (1) $a_1, a_2, a_3, a_4$ (2) $a_1, a_2, a_3, a_5$ (3) $a_1, a_2, a_3, a_6$ (4) $a_1, a_2, a_4, a_5$ (5) $a_1, a_2, a_5, a_6$ (6) $a_1, a_2, a_5, a_6$ (7) $a_1, a_3, a_4, a_5$ (8) $a_1, a_3, a_4, a_6$ (9) $a_1, a_3, a_5, a_6$ (10)

---

[9] G.W. Leibniz, Op. cit., p. 67
[10] Ebenhard Knobloch, Die Mathematischen Studien von G.W. Leibniz zur Kombinatorik, Franz Steiner Verlag, Wiesbaden 1973, p. 38; Michel Serres, Le systěme de Leibniz et ses modĕles mathématiques, P.U.F., 1968, p. 414.

$a_1$, $a_4$, $a_5$, $a_6$ (11) $a_2$, $a_3$, $a_4$, $a_5$ (12) $a_2$, $a_3$, $a_4$, $a_6$ (13) $a_2$, $a_3$, $a_5$, $a_6$ (14) $a_2$, $a_4$, $a_5$, $a_6$ (15) $a_3$, $a_4$, $a_5$, $a_6$.

(ii) 6 con5nations: ($a_1$, $a_2$, $a_3$, $a_4$ $a_5$), ($a_1$, $a_2$, $a_3$, $a_4$ $a_6$), ($a_1$, $a_2$, $a_3$, $a_5$ $a_6$), ($a_1$, $a_2$, $a_4$, $a_5$ $a_6$), ($a_1$, $a_3$, $a_4$, $a_5$ $a_6$), ($a_2$, $a_3$, $a_4$, $a_5$ $a_6$).

Hence, the total number of terms is

$$C_5^6 + C_5^6 \cdot C_1^5 = 6 + 6 \cdot 5 = 36.$$

(See Appendix II)

### III
### Considerations

This "mise au net" of Leibniz's insight uncovers a method to obtain a system of $k^2$ signs. We have shown that the notation employed gives "a local habitation and a name"[11] to the possible arrangements of the $k$ primitive elements. Therefore, combinations and semi-fractions might be envisaged as a syntax. A fact is worth emphasizing: in general, there are $k$ classes of $k$ sequences of symbols with identical signification. For instance the words ($a_1$, $a_2$, $a_3$, $a_4$ $a_5$), $\frac{1}{4}a_5$, $\frac{2}{4}a_4$, $\frac{4}{4}a_3$, $\frac{7}{4}a_2$ and $\frac{11}{4} a_1$, correspond to the same content or meaning[12]. From a semantic point view, these sequences are identical. Are we playing, then, with a set of empty symbols and a few rules to arrange them? The answer might be affirmative but, I think, the gist of Leibniz's thought points out a conceptual region where "the lack of content of a pure structure is compatible with the idea that something contentful is an instance of a pure structure"[13].

---

[11] W. Shakespeare, A Midsummer Night's Dream, V, i.
[12] Appendix II.
[13] This last phrase is our somewhat curtailed version of an argument due to Jeremy Butterfield (we refer the reader to: J. Butterfield "Our Mathematical Universe?" arxiv: 1406.4348, June 2014).

## Acknowledgments

This research has been partially supported by Fondecyt (Project 1120019).

Appendix I

65.2  Inveni omnes Termini primi ponantur in una classe, et designentur notis quibusdam; commodissimum erit numerari.

66.3  Inter Terminos primos ponantur non solum res, sed et modi sive respectus.

67.4  Cum omnes termini orti varient distantia a primis, prout ex pluribus Terminis primis componuntur, seu prout est exponens Complexionis, hinc tot classes faciendae, quot exponentes sunt, et in eandem classem conjiciendi termini, qui ex eodem numero primorum componentur.

68.5  Termini orti per com2nationem scribi aliter non poterunt, quam scribendo terminos primos, ex quibus componuntur; et quia termini primit signati sunt numeris, scribantur duo numeri duos terminos signantes.

69.6  At Termini orti per con3nationem aut alias majoris etiam exponentis Complexiones, seu Termini qui sunt in clase 3tia et sequentibus, singuli toties varie scribi possunt, quot habet complexiones simpliciter exponens ipsorum, spectatus non jam amplius ut exponens, sed ut numerus rerum. Habet hoc suum fundamentum in Usu IX, v.g. sunto termini primi his numeris signati 3.6.7.9.; sitque terminus ortus in classe tertia, seu per con3nationen compositus, nempe ex 3bus simplicibus 3.6.9. et sint in classe 2da combinations hae: (1)3.6. (2)3.7. (3)3.9. (4)6.7. (5)6.9. (6)7.9.; ajo terminum illum datum classis 3tiae scribi posse vel sic: 3.6.9. exprimendo omnes simplices; vel exprimendo unum simplicem, et loco caeterorum duorum simplicium scribendo com2nationem, v.g. sic: $\overline{2}^{1}\cdot 9$ vel: $\overline{2}^{3}\cdot 6$. vel sic: $\overline{2}^{5}\cdot 3$. Hae quasi-fractiones quid significent, mox dicetur. Quo autem classis a prima remotior, hoc variatio major. Semper enim termini classis antecedentis sunt quasi genera subalterna ad terminus quosdam variationis sequentis.

70.7  Quoties terminus ortus citatur extra suam classem, scribatur per modum fractionis, ut numerus superior seu numerator sit numerus loci in classe, inferior seu nominator, numerus classis. 8. Commodius est, in terminis ortis exponendis non omnes terminos primos, sed intermedios scribere, ob multitudinem, ex iis eos qui maxime cogitanti de re occurrunt. Verum omnes primos scriber est fundamentalius.

# Appendix II
# Language of 36 terms

| | | | | | |
|---|---|---|---|---|---|
| $a_1a_2a_3a_4a_5$ | $\frac{1}{4}a_5$ | $\frac{2}{4}a_4$ | $\frac{4}{4}a_3$ | $\frac{7}{4}a_2$ | $\frac{11}{4}a_1$ |
| $a_1a_2a_3a_4a_6$ | $\frac{1}{4}a_6$ | $\frac{3}{4}a_4$ | $\frac{5}{4}a_3$ | $\frac{8}{4}a_2$ | $\frac{12}{4}a_1$ |
| $a_1a_2a_3a_5a_6$ | $\frac{2}{4}a_6$ | $\frac{3}{4}a_5$ | $\frac{6}{4}a_3$ | $\frac{9}{4}a_2$ | $\frac{13}{4}a_1$ |
| $a_1a_2a_4a_5a_6$ | $\frac{4}{4}a_6$ | $\frac{5}{4}a_5$ | $\frac{6}{4}a_4$ | $\frac{10}{4}a_2$ | $\frac{14}{4}a_1$ |
| $a_1a_3a_4a_5a_6$ | $\frac{7}{4}a_6$ | $\frac{8}{4}a_5$ | $\frac{9}{4}a_4$ | $\frac{10}{4}a_3$ | $\frac{15}{4}a_1$ |
| $a_2a_3a_4a_5a_6$ | $\frac{11}{4}a_6$ | $\frac{12}{4}a_5$ | $\frac{13}{4}a_4$ | $\frac{14}{4}a_3$ | $\frac{15}{4}a_2$ |